\documentclass{amsart}
\usepackage{topproc}
\usepackage{amssymb,url,color}

\newtheorem*{ack}{Acknowledgements}

\normalsize

\author{Zhikuang Chen}
\address{Zhili College, Tsinghua University, }
\email{zk-chen22@mails.tsinghua.edu.cn}

\author{Huaiqing Zuo}
\address{Department of Mathematical Sciences; Tsinghua University}
\email{hqzuo@mail.tsinghua.edu.cn}
\usepackage[margin=1.2in]{geometry}

\title[Poles of the Real Archimedean Zeta Function]
{On the Poles of Real Archimedean Zeta Functions}

\begin{document}
	
	\linespread{1.3}
	\Large
	
	\begin{abstract}
		This paper studies the poles of the real Archimedean zeta function for a weighted homogeneous polynomial $f \in \mathbb{R}[x, y]$ with an isolated singularity at the origin. By applying a weighted blow-up, we derive the meromorphic continuation of $Z_{f,\varphi}$ to $\text{Re }s > -1$. This explicit expression yields a necessary and sufficient condition for a root $s \in (-1, 0)$ of the Bernstein–Sato polynomial $b_f(s)$ to be a pole of $Z_{f,\varphi}$. Unlike the complex case established by F. Loeser (1985), this condition may fail in certain obvious cases—such as when $f$ is odd or even in $x$, $y$, or $(x, y)$—so not all such roots necessarily become poles.

        
	\end{abstract}
	
	\maketitle
	
	\section{Introduction}
	Let $f \in k[x_1,...,x_n]$ where $k = \mathbb R$ or $\mathbb C$. One can associate to $f$ a collection of invariants. The Monodromy Conjecture predicts a profound relationship between the motivic zeta function of $f$ (an algebraic invariant) and its local monodromies (a topological invariant); see \cite{Loe98}. The Archimedean zeta function $Z_{f,\varphi}$ defined for both real and complex fields, is an important tool for studying topological invariants of singularities. In the complex case  $k = \mathbb C$, the exponentials of its poles yield precisely the eigenvalues of the local monodromy. We present the relevant results in Subsection \ref{BFArchMono}.

	In 1985, F. Loeser \cite{Loe85} conjectured that for a complex version of the zeta function $Z_{f,\varphi}(s)$, a point  $s_0\in(-1,0)$ can be a pole if and only if $s_0$ is a root of the Bernstein-Sato polynomial $b_f(s)$ (see Definition \ref{bf}). He  proved this result under the assumption that $f$ is weighted homogeneous with an isolated singularity at the origin. However, in 2024, Dougal Davis, András C. Lőrincz and Ruijie Yang \cite{RY24} used the theory of polarized complex Hodge modules to study the poles of the Archimedean zeta function. They pointed out that if $\alpha\in(0,1)$ satisfies the condition that $-\alpha$ is a simple root of $b_f(s)$ and $f^{-\alpha}\in\mathcal D_{\mathbb C^n}\cdot f^{-\alpha+1}$, then $-\alpha$ cannot be a pole of $Z_{f,\varphi}$. This shows that Loeser's conjecture does not hold in general. Moreover, the properties of poles of the Archimedean zeta function in the real case can be more complicated than those in the complex case. This paper therefore focuses on the poles of the real version of the Archimedean zeta function.
	
	We begin by reviewing results on the Bernstein-Sato polynomial. As such, this object was originally introduced to achieve the meromorphic continuation of the Archimedean zeta function \cite{Be}. 
	\begin{definition} \label{bf}(Bernstein-Sato polynomial) 
    Let $k$ be a field of characteristic zero and $f \in k[x_1, \dots, x_n]$ a polynomial. The theory of holonomic $D$-modules \cite{Cou95} guarantees the existence of a non-zero differential operator $$\mathcal{D} \in k[ s, x_1, \dots, x_n, \frac{\partial}{\partial x_1}, \dots, \frac{\partial}{\partial x_n} ]$$ and a polynomial $b(s) \in k[s] \setminus \{0\}$ satisfying the functional equation $\mathcal{D} \cdot f^{s+1} = b(s) \cdot f^s$. Therefore, the set
    $$I:=\{b(s)\in k[s]\mid \exists \;\mathcal D \; \text{such that}\;  \mathcal D\cdot f^{s+1}=b(s)\cdot f^s\} $$
    is a non-zero ideal of $k[s]$. Its monic generator is defined to be the Bernstein–Sato polynomial $b_f(s)$ of $f$.
	\end{definition}

When $k = \mathbb C$, a  fundamental result due to Kashiwara states that all roots of the Bernstein–Sato polynomial  $b_f(s)$ are negative rational numbers \cite{Rationality_of_Roots_of_B-Function_Kashiwara}. This rationality property was subsequently extended to the case $k = \mathbb R$ using \cite[Proposition 2.1]{BMS06}.

	\begin{definition} (The real Archimedean zeta function) 
    For a polynomial $$f\in\mathbb R[x_1,\cdots,x_n]$$ and a test function $\varphi\in C_c^\infty(\mathbb R^n)$ (i.e., a compactly supported, real-valued, smooth function on $\mathbb R^n$), the real Archimedean zeta function  is defined as  
		$$Z_{f,\varphi}(s):=\int_{\mathbb R^n}|f(\mathbf x)|^s\varphi(\mathbf x)\mathrm d\mathbf x,\text{Re }s>0.$$

        This integral is well-defined and, owing to the compact support of $\varphi$, defines a holomorphic function of $s$ in the half-plane $\text{Re }s>0$.
		\end{definition}
	
	Let $b_f(s)$ be the Bernstein-Sato polynomial of $f$ over $\mathbb R$, and choose a differential operator $\mathcal D$ satisfying  $\mathcal D\cdot f^{s+1}=b_f(s)\cdot f^s$. By substituting $-f$ for $f$ in this equation, we find that the same operator $\mathcal D$ and polynomial $b_f(s)$ satisfy the relation  $$(-\mathcal D)\cdot(-f)^{s+1}=b_f(s)\cdot(-f)^s.$$
	
	Since $\varphi$ has compact support in $\mathbb R^n$, we may apply integration by parts to obtain:
	
	$$Z_{f,\varphi}(s)=\int_{\mathbb R^n}|f|^s\varphi\mathrm d\mathbf x$$
	$$=\frac1{b_f(s)}\left(\int_{f>0}\mathcal D\cdot f^{s+1}\varphi\mathrm d\mathbf x+\int_{f<0}(-\mathcal D)\cdot(-f)^{s+1}\varphi\mathrm d\mathbf x\right)$$
	$$=\frac1{b_f(s)}\left(\int_{f>0}f^{s+1}\overline{\mathcal D}\cdot\varphi\mathrm d\mathbf x-\int_{f<0}(-f)^{s+1}\overline{\mathcal D}\cdot\varphi\mathrm d\mathbf x\right),$$
	
	where $\overline{\mathcal D}$ is the "dual operator" of $\mathcal D$, defined by the following relations:
	$$\overline{x_i}=x_i,\; \overline{\frac\partial{\partial x^i}}=-\frac\partial{\partial x^i},\; \overline{\mathcal D_1+\mathcal D_2}=\overline{\mathcal D_1}+\overline{\mathcal D_2},\; \overline{\mathcal D_1\cdot\mathcal D_2}=\overline{\mathcal D_2}\cdot\overline{\mathcal D_1}.$$
	
	Since the expression $\int_{f>0}f^{s+1}\overline{\mathcal D}\cdot\varphi\mathrm d\mathbf x-\int_{f<0}(-f)^{s+1}\overline{\mathcal D}\cdot\varphi\mathrm d\mathbf x$ is well-defined and holomorphic for $\text{Re }s>-1$, it provides a meromorphic continuation of $Z_{f,\varphi}(s)$ to this half-plane. Moreover, it's easy to see if $s$ is a pole of $Z_{f,\varphi}(s)$ for some $\varphi\in C_c^\infty(\mathbb R^n)$, then we must have $b_f(s)=0$.

\begin{remark}\label{rk1.3}
 An inductive argument shows that by incorporating the factors $b_f(s+k)$ for $k=0,...,m$, one obtains a meromorphic continuation of $Z_{f,\varphi}(s)$ to the region $\mathrm{Re}\, s > -m-1$. Now, suppose $s_0$ is a pole. The construction of the continuation implies that $b_f(s_0+k) = 0$ for some $k\in \{0,1,...,\lfloor s_0\rfloor+1\}$. This implies that $s_0=(s_0+k)-k$ is a translate of a root of $b_f(s)$, and since all roots of the Bernstein–Sato polynomial are negative rationals, the same must be true for all poles of $Z_{f,\varphi}(s)$.
\end{remark}
	
	\begin{remark}\label{rmk1} Similarly, the complex version of the Archimedean zeta function can be defined as follows (see \cite{Igu, RY24}): 
		$$Z_{f,\varphi}(s):=\int_{\mathbb C^n}|f(\mathbf x)|^{2s}\varphi(\mathbf x)\mathrm d\mu(\mathbf x),\; \text{Re }s>0,$$
		where $\mathrm d\mu(\mathbf x)$ denotes  the Lebesgue measure on $\mathbb C^n$, identified with $\mathbb R^{2n}$. This function admits a 
 meromorphic continuation to the half-plane $\text{Re }s>-1$, and in fact, to the entire complex plane $\mathbb C$, with its poles located at a discrete set of negative rational numbers; see \cite{Igu,BSZ25} for details.
	\end{remark}

	Conversely, one may ask the following: given a root  $s_0$ of $b_f(s)$ satisfying  $s_0 > -1$, does there exist a test function  $\varphi\in C_c^\infty(\mathbb R^n)$ such that $s_0$ is a pole of the meromorphic continuation of  $Z_{f,\varphi}(s)$ to the half-plane $\text{Re }s>-1$?

This paper establishes a necessary and sufficient condition for the existence of a test function $\varphi\in C_c^\infty(\mathbb R^2)$ such that a given root $s_0\in(-1,0)$ of the Bernstein–Sato polynomial $b_f(s)$ becomes a pole of the real Archimedean zeta function $Z_{f,\varphi}$,  in the case where $f\in\mathbb R[x,y]$ is a bivariate weighted homogeneous polynomial with an isolated singularity at the origin.

While Loeser \cite{Loe85} established the converse for weighted homogeneous polynomials with an isolated singularity in the complex setting, our result demonstrates a fundamental difference in the real case: this correspondence fails in general. Specifically, we show that even for such polynomials, not every root $s_0\in(-1,0)$ of $b_f(s)$ corresponds to a pole of $Z_{f,\varphi}$ for some $\varphi$. 

Furthermore, we apply our main theorem to three specific types of such polynomials—namely fewnomial singularities (\cite{YZ16}), $f=x^n+y^m,\; f=x^m+xy^n,\; \text{and}\; f=x^ny+xy^m$—and characterize some roots $s_0\in(-1,0)$ do not appear as poles in these cases (see the final section).

We now state the main result of this paper.
	
	\begin{theorem}\label{mt} Let $f\in\mathbb R[x,y]$ be a weighted homogeneous polynomial of type $(a,b; m)$ with an isolated singularity at the origin, meaning $f(\lambda^ax,\lambda^by)=\lambda^mf(x,y)$, for all $\lambda\in \mathbb R$. Let $s_0\in(-1,0)$ be a root of the Bernstein-Sato polynomial  $b_f(s)$. Then the following are equivalent:
		
		\begin{enumerate}
			\item There exists a test function $\varphi\in C_c^\infty(\mathbb R^2)$ such that $s_0$ is a pole of the (meromorphically continued to $\text{Re }s>-1$) real Archimedean zeta function
			$$Z_{f,\varphi}(s)=\int_{\mathbb R^2}|f(x,y)|^s\varphi(x,y)\mathrm dx\mathrm dy,\text{Re }s>0.$$
			
			\item There exist non-negative integers  $j,k\in\mathbb N$ such that $$-ms_0=(j+1)\cdot a+(k+1)\cdot b,$$ and the following integral does not vanish:
			$$\int_\mathbb R\left(|f(1,u)|^{s_0}+(-1)^j|f(-1,u)|^{s_0}\right)u^k\mathrm du\ne0.$$
		\end{enumerate}
	\end{theorem}
\begin{remark}
We believe there is  an analogous characterization in higher dimensions, where the poles of the zeta function correspond to specific arithmetic data of the weights, as established here.
    \end{remark}

\textcolor{black}{\begin{ack}
	 We thank  Quan Shi and Yongxin Xu for  helpful discussions.  H. Zuo was supported by BJNSF Grant 1252009 and NSFC Grant 12271280.
\end{ack}}
	
	\section{Preliminaries}
	\subsection{Bernstein-Sato polynomials, Archimedian zeta functions, and monodromy eigenvalues}\label{BFArchMono}
	
	let $f \in \mathbb C[x_1,...,x_n]\setminus \mathbb C$ be a non-constant polynomial. The local monodromy eigenvalues of $f$ are local embedded topological invariants arsing from the Milnor fibration structure of $f$ \cite{Milnor_Fibration_Monodromy_2}. These eigenvalues can be interpreted as those of the nearby cycles functor \cite{Bry86}. We denote by $S(f)$ the set of all local monodromy eigenvalues, which is a finite subset of the roots of unity.  We summerize the interrelations among  $b_f(s), Z_{f,\varphi}(s)$, and $S(f)$ in the following Theorem \ref{Triangle of three local invariant}.
	
	Conventionally, we regard $Z_f(s)=Z_{f,-}(s)$ as a distribution (or a generalized function) in the test function $\varphi$. We say that $s_0$ is  a pole of order at least $r$ if there exists a test function $\varphi \in C_c^{\infty}(\mathbb C^n)$ such that $s_0$ is a pole order at least $r$ for the meromorphic function $Z_{f,\varphi}(s)$.
    
	\begin{theorem}\label{Triangle of three local invariant}
		Notations are as above, then:

			\noindent(I) (Kashiwara-Malgrange \cite{Kashiwara_Malgrange_Theorem_Kashiwara,Kashiwara_Malgrange_Theorem_Malgrange}) If $\lambda \in S(f)$, then there exists a root $\alpha$ of $b_f(s)$ such that $e^{2\pi i \alpha} = \lambda$. Conversely, the exponential of each root of $b_f(s)$ is a monodromy eigenvalue.
			
			\noindent(II) (Barlet \cite{Bar84}) If $\lambda$ is a monodromy eigenvalue for $f$, that is, if $\lambda \in S(f)$, then there exists   a pole $\alpha$ of  $Z_{f}(s)$  such that $\lambda = e^{2\pi i \alpha}$.  In particular,		\begin{displaymath}
				\{ \lambda\in\mathbb C^*\mid  \lambda=e^{2\pi i \alpha}\text{ for some pole }\alpha\text{ of }Z_{f}(s)\} = S(f).
			\end{displaymath}
			
			\noindent(III) (Barlet) Suppose $\lambda\in S(f)$ is an eigenvalue of the monodromy on the cohomology $H^j$ of the Milnor fiber of $f$ 	
			at a point $x \in \{f = 0\}$, with a Jordan block of size  $m\ge 1$. Then:
			\begin{enumerate}
				\item (\cite{Bar84}) 
				There exists  a pole $\alpha$ of $Z_{f}(s)$ of order $\geq m$ such that $e^{2\pi i\alpha} = \lambda$. 
				\item (\cite{Bar86}) Write $\lambda=e^{2\pi i\alpha}$ with $\alpha\in (-1,0]\cap \mathbb{Q}$. If $j\ge 1$ then $\alpha-j$ is a pole of $Z_f(s)$ of order $\ge m$. 
				
				\item (\cite{Bar84b}) If $\lambda=1$ and $j\ge 1$, there exists  a negative integer pole  of $Z_{f}(s)$ of order $\geq m+1$.
			\end{enumerate}
	\end{theorem}
  
	\begin{remark}
		The last assertion of (2) follows from (1) and Remark \ref{rmk1}.
	\end{remark}

	\begin{remark}
		We point out all of the results above admit generalizations to an ordered collections of polynomials. However, $b_f(s),Z_{f,\varphi}(s),S(f)$ should be substituted by Bernstein-Sato ideals, multivariable Archimedean zeta functions, the support of Sabbah specialization complexes respectively, cf. \cite{Sab87,Alexander Module,Zero_Loci_of_Bernstein-Sato_Ideals,BS2,BSZ25}.
	\end{remark}
At the end of this subsection, we recall a result on the propagation of poles for $Z_f(s)$; this corresponds to the case $r=1$ in \cite[Lemma 2.13]{BSZ25}.
	
		\begin{lemma}\label{propagation_of_poles}
		Let $f\in \mathbb C[x_1,...,x_n]\setminus \mathbb C$. If $s_0$ is a pole of $Z_f(s)$ of order $\geq m$, then $s_0-1$ is again a pole of $Z_f(s)$ of order $\geq m$.
	\end{lemma}
	\begin{proof}
		Suppose $s_0$ is an order $m$ pole of $Z_{f,\omega}(s)$, where $\omega \in C_c^{\infty}(\mathbb C^n)$. Then $s_0-1$ is an order $m$ polar hyperplane of $Z_{f,\omega}(s+1) = Z_{f,\vert f\vert^2 \omega}(s)$, and hence is an order $\geq m$ pole of $Z_f$.
	\end{proof}
	\begin{remark}
		In the real setting, it appears less evident whether the same propagation holds, as $\vert f\vert$ is not necessarily smooth. However, the same proof suggests that if $s_0$ is a pole, then $s_0 - 2$ must also be one.
	\end{remark}
	
	\subsection{Roots of Bernstein-Sato polynomials in the weighted homogeneous case}

   The following standard result will be used.
	
	\begin{lemma}\label{2.6} If a polynomial $f\in\mathbb R[x_1,\cdots,x_n]$ is weighted homogeneous of type $(a_1,\cdots,a_n; m)$, i.e. $f(\lambda^{a_1}x_1,\cdots,\lambda^{a_n}x_n)=\lambda^mf(x_1,\cdots,x_n)$, with an isolated singularity at the origin (implying a finite-dimensional  Milnor algebra), then the following identity holds:
		$$b_f(s)=(s+1)\prod\limits_{[\mathbb R[x_1,\cdots,x_n]/(\partial f)]_t\ne0}\left(s+\frac{t+\sum\limits_{i=1}^na_i}m\right),$$
		where $[\mathbb R[x_1,\cdots,x_n]/(\partial f)]_t$ is the set of weighted homogeneous elements of weighted degree $t$ in the Milnor algebra $\mathbb R[x_1,\cdots,x_n]/(\partial f)=\mathbb R[x_1,\cdots,x_n]/(\partial_{x_1}f,\cdots,\partial_{x_n}f)$.
	\end{lemma}
	
	\begin{proof} See \cite{BGM86}.
	\end{proof}
	
	\begin{corollary}\label{2.7} Under the conditions of Lemma \ref{2.6}, the set of roots of  $b_f(s)$ in  $(-1,0)$ is given by 
		$$\{-\frac{\sum\limits_{i=1}^nk_ia_i}m\mid k_i\in\mathbb Z_+,\sum\limits_{i=1}^nk_ia_i<m\}.$$
	\end{corollary}
	
	\begin{proof} It is obvious that the ring $R=\mathbb R[x_1,\cdots,x_n]$ contains a nonzero weighted homogeneous element of weighted degree $t$  if and only if $t$ can be expressed as a non-negative integer linear combination of $a_1,\cdots,a_n$. 
    
    Moreover, since the ideal $(\partial_{x_1}f,\cdots,\partial_{x_n}f)$ is generated by weighted homogeneous polynomials of weighted degree at least $m-\max\{a_1,\cdots,a_n\}$, it follows that $[(\partial f)]_t=0$ for all $t<m-\sum\limits_{i=1}^n a_i$. This implies that if we set $s=-\frac{t+\sum\limits_{i=1}^na_i}m\in(-1,0)$, then $[R/(\partial f)]_t\ne0\Leftrightarrow t$ can be written as non-negative integer linear combination of $a_1,\cdots,a_n$. Equivalently, by Lemma \ref{2.6}, $s\in(-1,0)$ is a root of $b_f(s)$ if and only if $s=-\frac dm$, where $d$ can be expressed as a positive integer linear combination of $a_1,\cdots,a_n$.
	\end{proof}
	
	\section{Proof of Theorem \ref{mt}}
	By applying the transformation
	\begin{equation}
		\begin{cases}
			x=\pm t^a,
			\\y=t^bu,
		\end{cases}\; \text{where }t\ge0, \; u\in\mathbb R,
	\end{equation}
	
	we obtain
	
	\begin{equation}
		\begin{aligned}
			Z_{f,\varphi}(s)&=a\left(\int_0^\infty t^{a+b-1}\mathrm dt\int_{\mathbb R}|f(t^a,t^bu)|^s\varphi(t^a,t^bu)\mathrm du\right.
			\\&+\left.\int_0^\infty t^{a+b-1}\mathrm dt\int_{\mathbb R}|f(-t^a,t^bu)|^s\varphi(-t^a,t^bu)\mathrm du\right)
			\\&=a\left(\int_{\mathbb R}|f(1,u)|^s\mathrm du\int_0^\infty t^{ms+a+b-1}\varphi(t^a,t^bu)\mathrm dt\right.
			\\&+\left.\int_{\mathbb R}|f(-1,u)|^s\mathrm du\int_0^\infty t^{ms+a+b-1}\varphi(-t^a,t^bu)\mathrm dt\right)
			\\ &=\frac{(-1)^{m-a-b+1}a}{\prod\limits_{k=a+b}^m(ms+k)}\left(\int_{\mathbb R}|f(1,u)|^s\mathrm du\int_0^\infty t^{m(s+1)}\right.
			\\&\left.\frac{\mathrm d^{m-a-b+1}}{\mathrm dt^{m-a-b+1}}\varphi(t^a,t^bu)\mathrm dt+\int_{\mathbb R}|f(-1,u)|^s\mathrm du\int_0^\infty t^{m(s+1)}\right.
			\\&\left.\frac{\mathrm d^{m-a-b+1}}{\mathrm dt^{m-a-b+1}}\varphi(-t^a,t^bu)\mathrm dt\right).
		\end{aligned}
	\end{equation}
	\begin{remark} The transformation corresponds locally to a weighted blow-up (cf. \cite{Motivic_Zeta_Function_on_Q_Gorenstein_Varieties_Leon_Martin_Veys_Viu_Sos}). For weighted homogeneous polynomials with an isolated singularity at the origin, the weighted blow-up yields an embedded $\mathbb{Q}$-resolution (see \cite{Monodromy_Conjecture_for_Semi-Quasihomogeneous_Hypersurfaces_Blanco_Budur_vdV}).
	\end{remark}
	
	\begin{lemma} \label{3.2}
   Assuming that $f(x,y)$ is a weighted homogeneous polynomial of type $(a,b; m)$ with an isolated singularity at the origin, the polynomials  $f(1,u)$ and $f(-1,u)$, when viewed as functions of $u$, have no multiple roots in $\mathbb R$.
	\end{lemma}
	\begin{proof} 
    
     We firstly show the set $S=\{(x,y)\in\mathbb R^2\mid x\ne0,f(x,y)=f_y(x,y)=0\}$ is finite. Note that if $f$ is weighted homogeneous and has an isolated singularity at the origin, then we have $f=\frac1m\left(axf_x+byf_y\right)$. In this way, $f=f_y=0$ implies $x=0$ or $f_x=0$, so it suffices to show $\{(x,y)\in\mathbb R^2\mid f_x(x,y)=f_y(x,y)=f(x,y)=0\}=:\text{Sing}(f)$ is finite. This is immediate by the fact that the origin is an isolated singularity.
		
	Back to the lemma, if $f(1,u)$ has a multiple root $y\in\mathbb R$, then we have $f(1,u)=g(u)(u-y)^2$ for some $g\in\mathbb R[u]\Rightarrow f(1,y)=f_y(1,y)=0$. Since $f$ is a weighted homogeneous polynomial of type $(a,b; m)$ (and hence $f_y$ is a weighted homogeneous polynomial of type $(a,b;m-b)$), we have $f(\lambda^a,\lambda^by)=f_y(\lambda^a,\lambda^by)=0$ for any $\lambda\in\mathbb R^\times\Rightarrow S\supseteq\{(\lambda^a,\lambda^by)\mid\lambda\in\mathbb R^\times\}$, which contradicts $|S|<\infty$. Similarly for $f(-1,u)$.
    
	\end{proof}
	
	\begin{lemma}\label{3.3} For any given $a,b\in\mathbb Z^+$, There exist positive integers $c_{ijk}\in\mathbb Z^+,$ where $i,j,k\in\mathbb N$, depending only on $a$ and $b$, such that for all $N\in\mathbb N$ and all smooth functions $\varphi\in C^\infty(\mathbb R^2)$, the following identity holds:
		$$\frac{\mathrm d^N}{\mathrm dt^N}\varphi(\alpha t^a,t^bu)=\sum\limits_{\substack{0\le i+j\le N,k\ge0\\i\cdot a+j\cdot b=N+k}}c_{ijk}\alpha^iu^jt^k\frac{\partial^{i+j}}{\partial x^i\partial y^j}\varphi(\alpha t^a,t^bu),$$
		
	where $\alpha\in\mathbb R$ is any constant.
	\end{lemma}

	\begin{proof} We prove it by induction on $N$. To establish the base case when $N=0$, we take $c_{000}=1$, and the identity is clearly satisfied.
		
		Assume that the identity holds for $N=n$. For the case $N=n+1$, we  differentiate the inductive hypothesis with respect to $t$, obtaining:
\begin{align*}
    \frac{\mathrm d^{n+1}}{\mathrm dt^{n+1}}\varphi(\alpha t^a,t^bu)&=\sum\limits_{\substack{0\le i+j\le n,k\ge0\\i\cdot a+j\cdot b=n+k}}\left[kc_{ijk}\alpha^iu^jt^{k-1}\frac{\partial^{i+j}}{\partial x^i\partial y^j}\varphi(\alpha t^a,t^bu)\right.\\
    &\quad+ac_{ijk}\alpha^{i+1}u^jt^{k+a-1}\frac{\partial^{i+j+1}}{\partial x^{i+1}\partial y^j}\varphi(\alpha t^a,t^bu)\\
    &\quad+\left.bc_{ijk}\alpha^iu^{j+1}t^{k+b-1}\frac{\partial^{i+j+1}}{\partial x^i\partial y^{j+1}}\varphi(\alpha t^a,t^bu)\right].
\end{align*}
        
		
It is straightforward to verify that each term in the sum is of the form $$\alpha^pu^qt^r\frac{\partial^{p+q}}{\partial x^p\partial y^q}\varphi(\alpha t^p,t^qu)$$ and satisfies the condition $p\cdot a+q\cdot b=n+1+r$. 
        
From this expression, we can read off the recurrence relation for the coefficients $c_{ijk}$ with $0\le i+j\le n+1$ and $i\cdot a+j\cdot b=n+k+1$:
		$$c_{ijk}=kc_{ij(k+1)}+ac_{(i-1)j(k-a+1)}+bc_{i(j-1)(k-b+1)},$$
	with the convention that a coefficient is zero if any of its indices falls outside the admissible range.

    To show that $c_{ijk}\in\mathbb Z_+$, we consider two cases:
    
		If $i+j\le n$, then $c_{ij(k+1)}>0$ by induction.
        
        If $i+j=n+1$, then at least one of $c_{(i-1)j(k-a+1)}$ and $c_{i(j-1)(k-b+1)}$ is positive by the assumption of induction (in particular, we must have $k=i\cdot a+j\cdot b-n\ge(n+1)\cdot\min\{a,b\}-n\ge\min\{a,b\}$, hence at least one of $k-a+1$ and $k-b+1$ is positive). 
        
        In both cases,  $c_{ijk}$ is a positive integer. This completes the inductive step and the proof.

	\end{proof}

	\begin{lemma} \label{3.4}The function $\Phi_+(s)$, defined by 
		$$\Phi_+(s)=\int_{\mathbb R}|f(1,u)|^s\mathrm du\int_0^\infty t^{m(s+1)}\frac{\mathrm d^{m-a-b+1}}{\mathrm dt^{m-a-b+1}}\varphi(t^a,t^bu)\mathrm dt,$$
		
		is holomorphic in the half-plane $\text{Re }s>-1$.
	\end{lemma}
	
	\begin{proof} Let all distinct roots of $f(1,u)$ in $\mathbb R$ be $u_1,u_2,\cdots,u_N\in\mathbb R$. Take a sufficiently small $\varepsilon>0$ and a sufficiently large $L>0$. We can make decomposition:
		$$\Phi_+(s)=\sum\limits_{k=1}^N\int_{u_k-\varepsilon}^{u_k+\varepsilon}|f(1,u)|^s\mathrm du\int_0^\infty t^{m(s+1)}\frac{\mathrm d^{m-a-b+1}}{\mathrm dt^{m-a-b+1}}\varphi(t^a,t^bu)\mathrm dt$$
		$$+\int_{[-L,L]-\bigcup\limits_{k=1}^N(u_k-\varepsilon,u_k+\varepsilon)}|f(1,u)|^s\mathrm du\int_0^\infty t^{m(s+1)}\frac{\mathrm d^{m-a-b+1}}{\mathrm dt^{m-a-b+1}}\varphi(t^a,t^bu)\mathrm dt$$
		$$+\int_{|u|>L}|f(1,u)|^s\mathrm du\int_0^\infty t^{m(s+1)}\frac{\mathrm d^{m-a-b+1}}{\mathrm dt^{m-a-b+1}}\varphi(t^a,t^bu)\mathrm dt$$
		$$=:\sum\limits_{k=1}^N\Phi^1_k(s)+\Phi^2(s)+\Phi^3(s).$$
		
		We next assume that $\varphi$ is supported in $\overline{B(0,R)}$ and $\|\partial_x^p\partial_y^q\varphi\|_{L^\infty}\le M$ for every $p,q\in\mathbb N$ with $0\le p,q\le m-a-b+1$.
		
		Step 1: We show each $\Phi^1_k(s)$ is holomorphic in $\text{Re }s>-1$. By Lemma \ref{3.2}, we have $\frac1c|u-u_k|\le|f(1,u)|\le c|u-u_k|$ whenever $|u-u_k|\le\varepsilon$, for some $c>0$. Then when $\text{Re }s>-1$ we have:
		$$\int_{u_k-\varepsilon}^{u_k+\varepsilon}||f(1,u)|^s|\mathrm du\int_0^\infty\left|t^{m(s+1)}\frac{\mathrm d^{m-a-b+1}}{\mathrm dt^{m-a-b+1}}\varphi(t^a,t^bu)\right|\mathrm dt$$
		$$\le c^{|\text{Re }s|}R^{(m\text{Re }s+m+1)/a}\widetilde M\int_{u_k-\varepsilon}^{u_k+\varepsilon}|u-u_k|^{\text{Re }s}\mathrm du$$
		$$=2c^{|\text{Re }s|}R^{(m\text{Re }s+m+1)/a}\widetilde M\int_0^\varepsilon\xi^{\text{Re }s}\mathrm d\xi<\infty,$$
		
		where $\widetilde M:=\max\limits_{|u-u_k|\le\varepsilon,0\le t\le R^{1/a}}\left|\frac{\mathrm d^{m-a-b+1}}{\mathrm dt^{m-a-b+1}}\varphi(t^a,t^bu)\right|$.
		
		This shows $\Phi^1_k(s)$ is well-defined in $\text{Re }s>-1$.
		
		Also, if we define
		$$\Phi_\delta(s):=\int_{\delta<|u-u_k|<\varepsilon,t\le R^{1/a}}|f(1,u)|^st^{m(s+1)}\frac{\mathrm d^{m-a-b+1}}{\mathrm dt^{m-a-b+1}}\varphi(t^a,t^bu)\mathrm du\mathrm dt,$$
		then $\Phi_\delta(s)$ is holomorphic in $\text{Re }s>-1$ since the integral set is of finite measure and the integrated function is bounded and holomorphic with $s$, in the definition of $\Phi_\delta(s)$. Note that:
		$$|\Phi_k^1(s)-\Phi_\delta(s)|\le2c^{|\text{Re }s|}R^{(m\text{Re }s+m+1)/a}\widetilde M\int_0^\delta\xi^{\text{Re }s}\mathrm d\xi$$
		$$=2c^{|\text{Re }s|}R^{(m\text{Re }s+m+1)/a}\widetilde M\frac{\delta^{\text{Re }s+1}}{\text{Re }s+1}.$$
		
		This shows that for any $s_0$ with $\text{Re }s_0>-1$, there is a small neighborhood $U$ of $s_0$ such that $|\Phi_k^1(s)-\Phi_\delta(s)|$ uniformly tends to $0$ for all $s\in U$, as $\delta\to0^+$, which concludes that $\Phi_k^1(s)$ is holomorphic in $\text{Re }s>-1$.
		
		Step 2: We show $\Phi^2(s)$ is holomorphic in $\text{Re }s>-1$. In fact, since $|f(1,u)|$ has positive lower-bound and finite upper-bound in the compact set $[-L,L]-\bigcup\limits_{k=1}^N(u_k-\varepsilon,u_k+\varepsilon)$, the analyticity of $\Phi^2$ is clear.
		
		Step 3: It remains to prove $\Phi^3(s)$ is holomorphic in $\text{Re }s>-1$.
		
		Note that we can reduce the set of integral of $t$ to $[0,(R/|u|)^{1/b}]$ (since we only need to consider the set of $t$ with $(t^a,t^bu)\in B(0,R)$).
		
		Applying Lemma \ref{3.3}:
		$$\left|t^{m(s+1)}\frac{\mathrm d^{m-a-b+1}}{\mathrm dt^{m-a-b+1}}\varphi(t^a,t^bu)\right|\le Mt^{m(\text{Re }s+1)}\sum\limits_{\substack{0\le i+j\le m-a-b+1,k\ge0\\i\cdot a+j\cdot b=m-a-b+1+k}}c_{ijk}|u|^jt^k.$$
		
		Integrate on $[0,(R/|u|)^{1/b}]$:
		$$\int_0^{(R/|u|)^{1/b}}\left|t^{m(s+1)}\frac{\mathrm d^{m-a-b+1}}{\mathrm dt^{m-a-b+1}}\varphi(t^a,t^bu)\right|\mathrm dt$$
		$$\le M\sum\limits_{\substack{0\le i+j\le m-a-b+1,k\ge0\\i\cdot a+j\cdot b=m-a-b+1+k}}c_{ijk}|u|^j\int_0^{(R/|u|)^{1/b}}t^{m\text{Re }s+m+k}\mathrm dt$$
		$$\le M_0(s)\sum\limits_{\substack{0\le i+j\le m-a-b+1,k\ge0\\i\cdot a+j\cdot b=m-a-b+1+k}}|u|^{j-(m\text{Re }s+m+k+1)/b},$$
		
		where
		$$M_0(s):=\frac{MR^{(m\text{Re }s+m+k+1)/b}}{m\text{Re }s+m+k+1}\max\limits_{\substack{0\le i+j\le m-a-b+1,k\ge0\\i\cdot a+j\cdot b=m-a-b+1+k}}c_{ijk}.$$
		
		Let $D$ be the degree of $f(1,u)$, then we have $\frac1C|u|^D\le|f(1,u)|\le C|u|^D$ as $|u|\to\infty$, for some $C\in(0,\infty)$. It follows that:
		$$\int_{|u|>L}||f(1,u)|^s|\mathrm du\int_0^\infty\left|t^{m(s+1)}\frac{\mathrm d^{m-a-b+1}}{\mathrm dt^{m-a-b+1}}\varphi(t^a,t^bu)\right|\mathrm dt$$
		$$\le2C^{|\text{Re }s|}M_0(s)\sum\limits_{\substack{0\le i+j\le m-a-b+1,k\ge0\\i\cdot a+j\cdot b=m-a-b+1+k}}\int_L^\infty u^{j-(m\text{Re }s+m+k+1)/b+D\text{Re }s}\mathrm du.$$
		
		We hope to prove $j-\frac{m\text{Re }s+m+k+1}b+D\text{Re }s<-1$ for any $0\le i+j\le m-a-b+1,k\ge0$ with $i\cdot a+j\cdot b=m-a-b+1+k$, and in this way the above expression is well-defined and finite. Note that:
		$$j-\frac{m\text{Re }s+m+k+1}b+D\text{Re }s=\left(D-\frac mb\right)\text{Re }s+j-\frac{m+k+1}b$$
		$$=\left(D-\frac mb\right)\text{Re }s-\frac{(i+1)\cdot a}b-1=-\left(1+\frac{(i+1)\cdot a+(m-bD)\text{Re }s}b\right)\>\>\>\>(*).$$
		
		Since $f$ is weighted homogeneous of type $(a,b;m)$, the lowest power of $x$ should be $\frac{m-bD}a$ (note that $D$ is the largest power of $y$). Since the lowest power of $x$ can't be greater than $1$ (otherwise we have $x^2\mid f$, which contradicts Lemma \ref{3.2}), we have $m-bD=0$ or $a$, hence $(m-bD)\text{Re }s>-a$. This proves $(*)$ is smaller than $-1$, hence $\Phi^3$ is well-defined on $\text{Re }s>-1$.
		
		Similarly with Step 1, we can define:
		$$\Phi_T(s):=\int_{L<|u|<T,t\le R^{1/a}}|f(1,u)|^st^{m(s+1)}\frac{\mathrm d^{m-a-b+1}}{\mathrm dt^{m-a-b+1}}\varphi(t^a,t^bu)\mathrm du\mathrm dt.$$
		
		$\Phi_T(s)$ is holomorphic on $\text{Re }s>-1$, and for any $s_0$ with $\text{Re }s_0>-1$, there's a small neighborhood $U$ of $s_0$, such that $\Phi_T(s)$ uniformly tends to $\Phi^3(s)$ for all $s\in U$ as $T\to+\infty$. This concludes the proof.
	\end{proof}
	
	\begin{remark}\label{3.5} Similarly, the function
		$$\Phi_-(s):=\int_{\mathbb R}|f(-1,u)|^s\mathrm du\int_0^\infty t^{m(s+1)}\frac{\mathrm d^{m-a-b+1}}{\mathrm dt^{m-a-b+1}}\varphi(-t^a,t^bu)\mathrm dt$$
		
		is also holomorphic in the half-plane $\text{Re }s>-1$.
		
		In this way, (3.2) can be regard as another way to make a meromorphic continuation of $Z_{f,\varphi}(s)$ from $\text{Re }s>0$ to $\text{Re }s>-1$. That is,
		$$Z_{f,\varphi}(s)=\frac{(-1)^{m-a-b+1}a}{\prod\limits_{k=a+b}^m(ms+k)}(\Phi_+(s)+\Phi_-(s)).$$
		
		By the uniqueness of meromorphic continuation, this expression coincides with the continuation defined at the beginning. We now study the properties of $\Phi_+(s)+\Phi_-(s)$.
	\end{remark}
	
	\begin{lemma}\label{3.6}Fix $d\in\mathbb N$ and $i,j\in\{0,1,\cdots,d\}$. There exists a test function  $\varphi\in C_c^\infty(\mathbb R^2)$ whose derivatives at the origin satisfy
    $$\frac{\partial^{k+l}}{\partial x^k\partial y^l}\varphi(0,0)=\delta_{ik}\delta_{jl},\;\text{for all}\; k,l\in\{0,1,\cdots,d\},$$
    where $\delta$ denotes the Kronecker delta.
	\end{lemma}
	
	\begin{proof} For a fixed $d\in\mathbb N$ and $r\in\{0,1,\cdots,d\}$, we firstly construct some $\psi_{dr}\in C_c^\infty(\mathbb R)$ such that $\frac{\partial^k}{\partial x^k}\psi_{dr}(0)=\delta_{rk},\forall k\in\{0,\cdots,d\}$. This can be achieved since if we let $\rho(x)=\begin{cases}e^{-1/(1-x^2)}, &|x|\le 1\\0, &|x|>1\end{cases}$ and $E_i(x)=x^i\rho(x),i=0,1,\cdots,d$, then the matrix $E_{ij}=\frac{\partial^jE_i(0)}{\partial x^j}$ is upper-triangular and has non-zero diagonal elements, hence invertible. This shows we can represent $\psi_{dr}$ as a linear combination of $\{E_i\}_{i=0}^d$, namely $\psi_{dr}=\sum\limits_{i=0}^dc_iE_i$. By choosing $c_0,c_1,\cdots,c_d$ properly, our requirements can be satisfied.
		
		Back to the lemma, it suffices to let $\varphi(x,y)=\psi_{di}(x)\psi_{dj}(y)$.
	\end{proof}
	
	\subsection*{Proof of Theorem \ref{mt}.} By Lemma \ref{3.4} and Remark \ref{3.5}, if we fix $s_0=-\frac dm$, where $d\in\mathbb Z\cap[a+b,m)$, then the problem reduces to determining under what conditions there exists a function $\varphi\in C_c^\infty(\mathbb R^2)$ such that $\Phi_+(s_0)+\Phi_-(s_0)\ne0$.
	
	The first observation is that, applying integration by parts for $m-d$ times:
	$$\int_0^\infty t^{m(-\frac dm+1)}\frac{\mathrm d^{m-a-b+1}}{\mathrm dt^{m-a-b+1}}\varphi(\alpha t^a,t^bu)\mathrm dt$$
	$$=(-1)^{m-d}\int_0^\infty\frac{\mathrm d^{d-a-b+1}}{\mathrm dt^{d-a-b+1}}\varphi(\alpha t^a,t^bu)\mathrm dt$$
	$$=(-1)^{m-d-1}\left.\frac{\mathrm d^{d-a-b}}{\mathrm dt^{d-a-b}}\varphi(\alpha t^a,t^bu)\right|_{t=0},$$
	
	where $\alpha=\pm1,u\in\mathbb R$.
	
	Applying Lemma \ref{3.3}, when $s_0=-\frac dm$ we have:
	$$\Phi_+(s_0)+\Phi_-(s_0)$$
	$$=(-1)^{m-d-1}\int_{\mathbb R}\big(|f(1,u)|^{s_0}\sum\limits_{(i+1)\cdot a+(j+1)\cdot b=d}c_{ij0}u^j\frac{\partial^{i+j}}{\partial x^i\partial y^j}\varphi(0,0)$$
	$$+|f(-1,u)|^{s_0}\sum\limits_{(i+1)\cdot a+(j+1)\cdot b=d}c_{ij0}(-1)^iu^j\frac{\partial^{i+j}}{\partial x^i\partial y^j}\varphi(0,0)\big)\mathrm du$$
	$$=(-1)^{m-d-1}\sum\limits_{(i+1)\cdot a+(j+1)\cdot b=d}c_{ij0}\frac{\partial^{i+j}}{\partial x^i\partial y^j}\varphi(0,0)$$
	$$\cdot\int_{\mathbb R}\left(|f(1,u)|^{s_0}+(-1)^i|f(-1,u)|^{s_0}\right)u^j\mathrm du.$$
	
	If there exist some $i,j$ such that $(i+1)\cdot a+(j+1)\cdot b=d$ and $\int_{\mathbb R}(|f(1,u)|^{s_0}+(-1)^i|f(-1,u)|^{s_0})u^j\mathrm du\ne0$, then by Lemma \ref{3.6} we may find some $\varphi\in C_c^\infty(\mathbb R^2)$ such that $\frac{\partial^{i+j}}{\partial x^i\partial y^j}\varphi(0,0)=1$ and $\frac{\partial^{k+l}}{\partial x^k\partial y^l}\varphi(0,0)=0,\forall(k,l)\in\{0,\cdots,d\}^2-\{(i,j)\}$ are satisfied. Since $c_{ij0}>0$, the chosen $\varphi$ satisfies our requirements. Conversely, if for any $i,j$ with $(i+1)\cdot a+(j+1)\cdot b=d$ one has $\int_{\mathbb R}(|f(1,u)|^{s_0}+(-1)^i|f(-1,u)|^{s_0})u^j\mathrm du=0$, then we must have $\Phi_+(s_0)+\Phi_-(s_0)=0$ for any $\varphi\in C_c^\infty(\mathbb R^2)$, i.e. $s_0$ cannot be a pole of $Z_{f,\varphi}(s)$.
	
	\section{A comparison with monodromy eigenvalues}
    Recall from Theorem \ref{Triangle of three local invariant} that  the set of monodromy eigenvalues is characterized in two equivalent ways:
\begin{align*}
        S(f)&=\{e^{2\pi i\alpha}\mid\alpha\text{ is a root of }b_f(s)\}\\
        &=\{e^{2\pi i\alpha}\mid \exists \varphi\in C_c^\infty(\mathbb R^n), \text{ s.t. }\alpha\text{ is a pole of } Z_{f,\varphi}(s)\}.
    \end{align*}

    For a given  $\lambda\in S(f)$, we define the following two values:
    $$\alpha_1(\lambda)=\max\{\alpha<0\mid e^{2\pi i\alpha}=\lambda,\; \alpha\text{ is a root of }b_f(s)\},$$
    $$\alpha_2(\lambda)=\max\{\alpha<0\mid e^{2\pi i\alpha}=\lambda, \; \alpha\text{ is a pole of }Z_{f,\varphi}(s)\text{ for some }\varphi\in C_c^\infty(\mathbb R^n)\}.$$

    It is clear that both $\alpha_1(\lambda)$ and $\alpha_2(\lambda)$ are well-defined. Moreover, by Remark \ref{rk1.3} one readily observes $\alpha_2(\lambda)\le\alpha_1(\lambda)$. Under the settings of Theorem 1.5, we distinguish the following cases for $\lambda\ne1$:

\begin{enumerate}
\item $\alpha_1(\lambda) = \alpha_2(\lambda) \in (-1, 0)$: This case occurs when $\alpha = \alpha_1(\lambda)$ is both a root of $b_f(s)$ and a pole of $Z_{f,\varphi}(s)$ for some $\varphi\in C_c^\infty(\mathbb R^2)$.

\item $\alpha_2(\lambda) \in (-\infty, -1)$ and $\alpha_1(\lambda) \in (-1, 0)$ with $\alpha_2(\lambda) < \alpha_1(\lambda)$: Here, $\alpha = \alpha_1(\lambda)$ is a root of the Bernstein–Sato polynomial $b_f(s)$ but not a pole of $Z_{f,\varphi}(s)$ for any $\varphi\in C_c^\infty(\mathbb R^2)$. Moreover, there exists a minimal $j \in \mathbb{N}_+$ and $\varphi\in C_c^\infty(\mathbb R^2)$ such that $\alpha - j = \alpha_2(\lambda)$ is a pole of $Z_{f,\varphi}(s)$.

\item $\alpha_1(\lambda) < -1$: In this case, the monodromy eigenvalue $\lambda$ cannot be obtained as $e^{2\pi i \alpha}$ for any root $\alpha$ of $b_f(s)$ in $(-1, 0)$. Note that roots of the Bernstein–Sato polynomial in $(-\infty, -1)$ may not admit as explicit a description as those in $(-1, 0)$, warranting further investigation.
\end{enumerate}



    
	
	\section{Examples}
	
	Although the conditions in Theorem \ref{mt} appear non-trivial, there are simple cases where  some  roots of $b_f(s)$ in $(-1,0)$ cannot be poles of $Z_{f,\varphi}(s)$. Indeed, it is straightforward  to construct  examples when  $f$ satisfies any of the following:
	
	(1) $f(x,y)=f(-x,y)$;
	
	(2) $f(x,y)=-f(-x,y)$;
	
	(3) $f(x,y)=f(x,-y)$;
	
	(4) $f(x,y)=-f(x,-y)$;
	
	(5) $f(x,y)=f(-x,-y)$;
	
	(6) $f(x,y)=-f(-x,-y)$.
	
Assume that condition (5) is satisfied. Then we have:
\begin{align*}
    \int_{\mathbb R}&(|f(1,u)|^{s_0}+(-1)^i|f(-1,u)|^{s_0})u^j\mathrm du\\
    &=\int_0^\infty u^j(|f(1,u)|^{s_0}+(-1)^i|f(-1,u)|^{s_0})\mathrm du\\
    &+(-1)^j\int_0^\infty u^j(|f(-1,u)|^{s_0}+(-1)^i|f(1,u)|^{s_0})\mathrm du\\
    &=\int_0^\infty u^j((1+(-1)^{i+j})|f(1,u)|^{s_0}+((-1)^i+(-1)^j)|f(-1,u)|^{s_0})\mathrm du.
\end{align*}

Now, using the notations from  Theorem \ref{mt}, suppose  $0<d<m$.  If for every representation of $d$ in the form $d=(i+1)\cdot a+(j+1)\cdot b$ with $i,j\in\mathbb N$, the sum $i+j$ is odd (this holds, for example, when $d=a+b+\min\{a,b\}$), then the expression above vanishes identically. This implies that  $s_0=-\frac dm$ cannot be a pole of $Z_{f,\varphi}(s)$.

	\begin{example} \label{5.1} $f=x^n+y^m$.

    The polynomial $f$ is weighted homogeneous of type $(m,n; nm)$. By Corollary \ref{2.7}, the set of  roots of its Bernstein-Sato polynomial $b_f(s)$ in the interval $(-1,0)$ is given by $\{-\frac{mj+nk}{nm}\mid j,k\in\mathbb Z_+,mj+nk<nm\}$.

    Case 1: $m$  odd,  $n$  even

    In this case, the symmetry $f(x,y)=f(-x,y)$ holds. It follows that 

    $$\int_\mathbb R(|f(1,u)|^{s_0}+(-1)^i|f(-1,u)|^{s_0})u^j\mathrm du=\int_\mathbb R(1+(-1)^i)|f(1,u)|^{s_0}u^j\mathrm du.$$

    Suppose $s_0$ is such that for every representation
    $$-nms_0=(i+1)m+(j+1)n,\;i,j\in\mathbb N,$$ the exponent $i$ is odd. Then $s_0$ cannot be a pole of $Z_{f,\varphi}$. In particular, if $s_0\in(-1,0)$ is a root of $b_f(s)$ and $-nms_0$ is even, then $s_0$ is not a pole of $Z_{f,\varphi}(s)$.

Example: let $f=x^4+y^3$. Then $s_0=-\frac56=-\frac{3\cdot2+4\cdot1}{12}$ is a root of Bernstein-Sato polynomial, but not a pole of $Z_{f,\varphi}(s)$.

  Case 2: $m$ even,  $n$  odd
  
Here we have the symmetry  $f(x,y)=f(x,-y)$.  A similar argument shows that if in every representation  $$-nms_0=(i+1)m+(j+1)n,\; i,j\in\mathbb N,$$
the exponent $j$ is odd, then $s_0$  is not a pole of $Z_{f,\varphi}(s)$. Again, this occurs when $-nms_0$ is even.

Case 3: $m$ and $n$ both odd or both even

  In this situation, the symmetry is either $f(x,y)=f(-x,-y)$ or $f(x,y)=-f(-x,-y)$. As discussed at the beginning of this section, the value
  $$s_0=-\frac{n+m+\min\{n,m\}}{nm}$$
  is a root of Bernstein-Sato polynomial (provided $n+m+\min\{n,m\}<nm$), but not a pole of $Z_{f,\varphi}(s)$.

Example: Let $f=x^3+y^5$. Then $s_0=-\frac{11}{15}=-\frac{5+2\cdot3}{15}$ is not a pole of $Z_{f,\varphi}(s)$.
\end{example}
	
	\begin{example} $f=x^m+xy^n$.

   The polynomial $f$ is weighted homogeneous of type $(n,m-1; mn)$ and the roots of $b_f$ in $(-1,0)$ are given by $\{-\frac{nj+(m-1)k}{nm}\mid j,k\in\mathbb Z_+,nj+(m-1)k<nm\}$.

        Case 1: Both $m,n$ are odd

        In this case, the function satisfies the  symmetry $f(x,y)=-f(-x,y)$. By  analogy with Example \ref{5.1}, if in every representation of $d$ as

        $$d=(i+1)n+(j+1)(m-1),\>i,j\in\mathbb N,$$

        the index $i$ is odd, then $-\frac d{nm}$ cannot be a pole of $Z_{f,\varphi}(s)$. In particular, if $d$ is even, then $-\frac d{nm}$ is not a pole. This is because $m-1$ is even and $n$ is odd, which forces $i+1$ to be even (i.e., $i$ is odd) in any such representation.

        Case 2: Both $m,n$ are even

        Here, $f(x,y)=f(x,-y)$ holds. Analogously, $-\frac d{nm}$ cannot be a pole  of $Z_{f,\varphi}(s)$ if, in every representation

        $$d=(i+1)n+(j+1)(m-1),\>i,j\in\mathbb N,$$

        the index $j$ is odd. In particular, this condition is always satisfied when  $d$ is even.

        Case 3: $n+m$ is odd

        In this situation, $f$ satisfies either $f(x,y)=f(-x,-y)$ or $f(x,y)=-f(-x,-y)$. For a given $d$ represented as $d=(i+1)n+(j+1)(m-1)$, we have

        $$\int_{\mathbb R}(|f(1,u)|^{-d/nm}+(-1)^i|f(-1,u)|^{-d/nm})u^j\mathrm du$$
        $$=\int_0^\infty u^j((1+(-1)^{i+j})|f(1,u)|^{-d/nm}+((-1)^i+(-1)^j)|f(-1,u)|^{-d/nm})\mathrm du$$

        If, in every such representation with $i,j\in\mathbb N$, the sum $i+j$ is odd, then the integrand vanishes identically.  Consequently, when $n+m$ is odd and $d=n+m-1+\min\{n,m-1\}$, the value $-\frac d{nm}$ cannot be a pole of $Z_{f,\varphi}(s)$.

	\end{example}
	
	\begin{example} $f=x^ny+xy^m$, where $n, m>1$.
		
		The argument in this case is quite similar to that in Example \ref{5.1}. The polynomial $f$ is weighted homogeneous of type $(m-1,n-1; nm-1)$, and the set of all roots of its Bernstein-Sato polynomial $b_f(s)$ in $(-1,0)$ is  $$\{-\frac{(m-1)j+(n-1)k}{nm-1}\mid j,k\in\mathbb Z_+,(m-1)j+(n-1)k<nm-1\}.$$

Case 1: $m$ odd, $n$ even

In this case, $f(x,y)=-f(x,-y)$.  By  analogy with  Example \ref{5.1}, if for every representation $$d=(i+1)(m-1)+(j+1)(n-1),\; i,j\in\mathbb N,$$
the index $j$ is odd, then  $-\frac d{nm-1}$ cannot be a pole of $Z_{f,\varphi}(s)$. In particular, if $d$ is even, then $-\frac d{nm-1}$ is not a pole. This follows because $m-1$ is even and $n-1$ is odd, which forces $j+1$ to be even (i.e., $j$ is odd) in any such representation.

Case 2: $n$ odd, $m$ even

Here $f(x,y)=-f(-x,y)$. Similarly, $-\frac d{nm-1}$ cannot be a pole when in every representation 
 $$d=(i+1)(m-1)+(j+1)(n-1),\; i,j\in\mathbb N,$$ the index $i$ is odd. This condition is automatically satisfied when $d$ is even.

 In summary, when  $n+m$ is odd and $d\in(1,nm-1)$ is even, $-\frac d{nm-1}$ cannot be a pole of  $Z_{f,\varphi}(s)$.

Case 3: $n, m$ both odd or both even

In this situation, $f$ satisfies either $f(x,y)=f(-x,-y)$ or $f(x,y)=-f(-x,-y)$. For $d=(i+1)(m-1)+(j+1)(n-1)$, we have:
\begin{align*}
    &\int_{\mathbb R}(|f(1,u)|^{-d/(nm-1)}+(-1)^i|f(-1,u)|^{-d/(nm-1)})u^j\mathrm du\\
    =&\int_0^\infty u^j((1+(-1)^{i+j})|f(1,u)|^{-d/(nm-1)}+((-1)^i+(-1)^j)|f(-1,u)|^{-d/(nm-1)})\mathrm du.
\end{align*}
If for every representation $d=(i+1)(m-1)+(j+1)(n-1),\; i,j\in\mathbb N$, the sum $i+j$ is odd, then the expression above vanishes identically. Consequently, when  $n+m$ is even and $d=n+m+\min\{n,m\}-3$, the value
$-\frac d{nm-1}$ cannot be a pole of $Z_{f,\varphi}(s)$.

	\end{example}

\end{document}